\documentclass[reqno, 11pt]{amsart}  

\usepackage{color}    
\usepackage{mathtools} 
\usepackage{amssymb}

\makeatletter
\def\ckech{\mathaccent"\accentclass@014}
\def\tah{\mathaccent"\accentclass@05E}
\makeatother

  \renewcommand\check{\bm\ckech} 
\renewcommand\hat{\bm\tah}

\newcommand\mysection[1]{%
\section{#1}\setcounter{equation}{0}}

\makeindex

\usepackage{bm}

\newtheorem{theorem}{Theorem}[section]
\newtheorem{lemma}[theorem]{Lemma}
\newtheorem{corollary}[theorem]{Corollary}

\theoremstyle{definition}

\theoremstyle{remark}

\newtheorem{remark}[theorem]{Remark}
 
\newcommand\loc{\textnormal{loc}}
\newcommand{\shharp}{=\!\!\!\!\|}  
\newcommand{\vsharp}{\asymp\kern -.5em\|}

 \makeatletter 
 \def\dashint{\operatorname{\,\,\,\mathclap{\!\int}\! 
 \!\text{\bf--}\!\!}}  
 \makeatother

\newcommand\sfp{{\sf p}}

\newcommand\bC{\mathbb{C}}
 
\newcommand\bM{\mathbb{M}}
\newcommand\bR{\mathbb{R}}

\newcommand\bS{\mathbb{S}}

\newcommand\frB{\mathfrak{B}}

\newcommand\cL{\mathcal{L}}

\newcommand\E{{\sf{E}}}
\renewcommand\L{{\sf{L}}}  {{}}
\newcommand\W{{\sf{W}}}

\renewcommand\){{\rm)}}

\def\+){\tmspace+\thinmuskip{.05em}\)}

\def\dashnorm{\,\,\text{\bf--}\kern-.5em\|}

\newcommand{\osc}{\operatornamewithlimits{osc}}

\renewcommand{\eqref}[1]{\text{(\ref{#1})}}

\begin{document}

\title[]{On solvability
of parabolic equations with singular coefficients in  odd mixed-norm
Morrey-Sobolev spaces}
\author[]{N.V. Krylov}
\address{School of Mathematics, University of Minnesota, Minneapolis, MN, 55455}
\email{nkrylov@umn.edu}

\renewcommand{\subjclassname}{\textup{2010} 
Mathematics Subject Classification}
 


\begin{abstract}
 We prove an existence and uniqueness
theorem for second-order parabolic
equations in the whole space  with constant zeroth-order coefficient in mixed-norm Morrey-Sobolev
spaces. The main coefficient $a$ is assumed to be measurable in $t$ and BMO in $x$ and the first-order coefficients $b$ are in
an appropriate mixed-norm Morrey  classes (thus admitting rather rough singularities). The mixed-norm Morrey-Sobolev spaces are ``odd'' in the sense
that the interior integration in the formula
defining the norm is performed with respect to
$t$ and not to $x$ as is customary.
 
\end{abstract}

\maketitle

\mysection{Introduction}

Let $\bR^{d}$ be a Euclidean space of points
$x=(x^{1},...,x^{d})$, $\bR^{d+1}=\{(t,x):
t\in\bR,x\in \bR^{d}\}$, $\bR^{d+1}_{t}=(t,\infty)\times\bR^{d}$.
In \cite{Kr_25} the author proved a theorem
about unique solvability in mixed-norm Morrey-Sobolev
spaces for second-order parabolic equations
in the whole space $\bR^{d+1}$. The main coefficient 
$a$ was assumed to be measurable in $t$ and BMO in $x$ and the lower-order coefficients $b$ and $c$  were in appropriate mixed-norm Morrey  classes (thus admitting rather rough singularities). Mixed-norm spaces in
\cite{Kr_25} are based on the norm
\begin{equation}
                           \label{8.22.1}
\|f\|_{L_{p,q}}:=\Big(\int_{\bR}\Big(\int_{\bR^{d}}|f|^{p}\,dx
\Big)^{q/p}\,dt\Big)^{1/q}.
\end{equation}
In contrast in this paper we use
\begin{equation}
                           \label{8.22.2}
\|f\|_{\L_{q,p}}:=\Big(\int_{\bR^{d}}\Big(\int_{\bR }|f|^{q}\,dt
\Big)^{p/q}\,dx\Big)^{1/p}.
\end{equation}
This is done for two reasons. First,
one can find in the literature
not completely justified claims that
much of the theory
using \eqref{8.22.2} is quite parallel to that
using \eqref{8.22.1}. Here,   indeed,
our main result, proved with all necessary details, looks almost identical to
the main result of \cite{Kr_25}. The second reason is that we need the result 
associated with exactly \eqref{8.22.2} in our investigation of It\^o stochastic equations.
Remark \ref{remark 5.10,3} sheds some light
to that effect.

In addition to these differences the exposition
here follows quite a different path. If in 
\cite{Kr_25}   obtaining main a priory estimates 
in {\em Morrey-Sobolev\/} spaces was tied to
the solvability of model equations in cylinders, which is a very nontrivial
fact, here we use only equations in the whole space $\bR^{d+1}$. This streamlines the matter and shows the path
for generalization of our results
on the solvability in {\em Morrey-Sobolev\/} spaces
to higher-order parabolic equations.
We leave these generalizations to the interested reader. In addition, we only
consider the operators with constant zeroth-order
coefficient. Again we leave  to the reader making it to belong to appropriate Morrey classes.

To compare our results with existing literature
observe that it looks like this article is the first one
 treating from the beginning to the end equations in mixed-norm Morrey spaces
based on \eqref{8.22.2}, in addition, to singular first-order coefficients
(this article is based on an earlier preprint
\cite{Kr_25_2} where $\beta\in(1,2)$).
Then,
in our setting we can allow $|b|=c/|x|$
with sufficiently small $c$.
This function  is independent of time
and it does not fit in \cite{FHS_17}, \cite{CLV_96},
\cite{MT_13}, which
are about the {\em elliptic\/} equations with Morrey coefficients and solutions in Morrey class
in \cite{FHS_17} or in Sobolev class in \cite{CLV_96}
and \cite{MT_13}. In the case of elliptic equation
 one can find plenty of information 
on generalized versions of Morrey spaces
and the Dirichlet problem  in the recent article 
\cite{GO_22} (with no lower-order terms 
in the main text) and in the references therein. 

The literature on the parabolic equations in Morrey spaces (no mixed-norms)
is not as rich as in the case of elliptic equations,
although it is worth drawing the reader's attention to
\cite{Li_03}, \cite{So_11}, \cite{Zh_19}, \cite{ZJZ_09}
and the references therein.
 However, it is worth
noting that equations in
the ``unusual'' mixed normed Sobolev (not Morrey-Sobolev)
spaces were earlier considered
in \cite{KN_14} (see also the references therein).

We use the following basic notation. Set $|G|$ 
to be the Lebesgue 
measure of $G\subset \bR^{d+1}$,
$$
\dashint_{G}f\,dxdt=\frac{1}{|G|}\int_{G}f\,dxdt,
\quad
\dashnorm f\| _{L_{p}(G)}=\Big(\dashint_{G}|f(t,x)|^{p}\,dxdt
\Big)^{1/p}.
$$
Similar notation is used for the functions 
of $x$ and $t$ only.

For $p,q\in(1,\infty)$ we 
introduce the spaces
$\L_{q,p} $, $\L_{q,p}(\Gamma)$,
as the sets of functions with finite norms
\eqref{8.22.2},
$$
\| f\|_{\L_{q,p }(\Gamma) }
=\| fI_{\Gamma}\|_{\L_{q,p } }.
$$
Define $\rho(t,x)=|x|+\sqrt{|t|}$,
$$
B_{R}=\{x:|x|<R\},\quad B_{R}(x)=x+B_{R},
$$
$$
C_{T,R}=[0,T)\times B_{R}, 
\quad C_{R}=
C_{R^{2},R},
$$
$$
C_{T,R}(t,x)=(t,x)+C_{T,R},
 \quad C_{R}(t,x)=C_{R^{2},R}(t,x),
$$
and let $\bC_{R}$ be the collection of $C_{R}(t,x)$ and $\bC=\cup_{R}\bC_{R}$.

Set
$$
\dashnorm f\|_{\L_{q,p}(C_{R}(t,x))}^{p}=
\dashint_{B_{R}(x)}\Big(
\dashint_{(t,t+R^{2})}|f(s,y)|^{q}\,ds\Big)^{p/q}\,dy.
$$
We use the notation $D_{i}=\partial/\partial x^{i}$, $D_{ij}=D_{i}D_{j}$,  
$Du=(D_{i}u)$, $D^{2}u=(D_{ij}u)$, $D^{n}u$
is the collection of all derivatives of $u$
of order $n$ with respect to $x$, 
$\partial_{t}u=\partial u/\partial t$.
All derivatives are assumed to be Sobolev
derivatives (if they exist). Denote
$C_{0}^{\infty}=C_{0}^{\infty}(\bR^{d+1})$.

Finally, fix $\delta\in(0,1]$ throughout the
article and set
$\bS_{\delta}$ to be the set of $d\times d$
symmetric matrices whose eigenvalues
are in $[\delta,\delta^{-1}]$.

\mysection{Main result}

Take $p,q\in(1,\infty)$ and
   $\beta\geq 0$ and introduce the
Morrey space $\E_{q,p,\beta} $ 
as the set of $g\in  \L_{q,p,\loc}$ 
such that  
\begin{equation}
                             \label{8.11.020}
\|g\|_{\E_{q,p,\beta} }:=
\sup_{\rho\leq 1,C\in\bC_{\rho}}\rho^{\beta}
\dashnorm g  \|_{ \L_{q,p}(C)} <\infty .
\end{equation} 
Observe that  $\E_{q,p,\beta}=\E_{q,p,2/q+d/p}$
for $\beta\geq 2/q+d/p$.
Define the Morrey-Sobolev spaces by
$$
\E^{1,2}_{q,p,\beta} =\{u:u,Du,D^{2}u,
\partial_{t}u\in \E_{q,p,\beta} \},
$$
where $Du,D^{2}u,
\partial_{t}u$ are Sobolev derivatives,
and 
provide $\E^{1,2}_{q,p,\beta} $ with an obvious norm.

Suppose that on $\bR^{d+1}$ we are given
an $\bS_{\delta}$-valued function $a(t,x)$
and an $\bR^{d}$-valued function $b(t,x)$.
 
For $\rho>0$  introduce  
\begin{equation}
                         \label{6.3.1}
 a^{\shharp}_{\rho} = \sup
_{r\leq\rho}
 \sup
_{ C\in\bC_{r} }\dashint_{C}|a(t,x)-\tilde a_{C}(t)|\,dxdt  ,\quad \tilde a_{C}(t)=\dashint_{C}a(t,x)\,dxds
\end{equation} 
(note $t$ and $ds$),
\begin{equation}
                           \label{3.14.2}
\bar b_{q ,p ,\rho }=\sup_{r\leq\rho }r
\sup_{C\in \bC_{r}} 
\dashnorm b \|_{\L_{q ,p  }(C)}.
\end{equation}

Define
$$
\cL u=\partial_{t}u+a^{ij}D_{ij}u+b^{i}u.
$$

Fix  $p ,q ,\beta    $ such that
\begin{equation}
                        \label{3.21.01}
p ,q  \in(1,\infty),\quad 1<\beta\leq\frac{d}{p }+\frac{2}{q }.
\end{equation}

Also fix $\rho_{a},\rho_{b}\in(0,1]$.
Here is our main result.

 \begin{theorem}
                     \label{theorem 5.8,20}
There 
exist
$$
\hat a=\hat a(d,\delta,p,q,\beta)>0,\quad\hat b=\hat b(d,\delta,p,q,\beta,\rho_{a} )>0,
$$ 
$$\lambda_{0}= \hat \lambda(d,\delta,p,q,\beta,\rho_{a} ) >0,\quad
N=N(d,\delta,p,q,\beta,\rho_{a}  ),
$$
such that, if $\lambda\geq \lambda_{0}\rho_{b}^{-2}$ and
$$
a^{\shharp}_{\rho_{a}}\leq 
\hat  a,\quad\bar b_{q\beta,p\beta,\rho_{b}}\leq \hat b,
$$
 then for any $u\in \E^{1,2}_{q,p,\beta}$ 
\begin{equation}
                        \label{5.10.2}
\|\lambda u,\sqrt\lambda Du, \partial_{t}u, D^{2}u\|_{  \E_{q,p,\beta}}
\leq N \rho_{b}^{-\alpha}\|f\|_{  \E_{q,p,\beta}}, 
\end{equation}
where   
$$
f=\cL u -\lambda u,\quad \alpha=
d+2+\beta-\frac{d}{p}-\frac{2}{q}.
$$
Furthermore, for any
$f\in \E_{q,p,\beta}$ and $\lambda\geq\lambda_{0}\rho_{b}^{-2}$ there exists a unique
$u\in \E^{1,2}_{q,p,\beta}$ such that in
$\bR^{d+1}$
$$
\cL u -\lambda u=f.
$$
 
\end{theorem}

This theorem is proved in Section \ref{subsection 8.22.1} after we develop
the necessary tools in Sections \ref{section
8.23.1} and \ref{section 4.29.1}.

\begin{remark}
                      \label{remark 5.10,3}
This theorem allows $b$ to be bounded. Also
the function 
$$
g(t,x):=\frac{1}{|x|+\sqrt{|t|}}I_{|x|\le	 1, |t| \leq1}
$$
  belongs, for instance, to
$\E_{p,p,1}$ as long as $1<p< d+2  $,
so that, if $|b|=cg$ with sufficiently small
$c$, $b$ is a valid choice for Theorem \ref{theorem 5.8,20} for $\rho_{b}=1$ and appropriate $\beta$. 

Indeed, if $|x|,\sqrt{|t|}\leq
2\rho$, then $C_{\rho}(t,x)\subset
(-5\rho^{2},5\rho^{2})\times B_{3\rho}$ and 
$$
\rho^{d+2}\dashnorm g\|^{p}_{\L_{p,p}( C_{\rho}(t,x))}=N\|g\|^{p}_{\L_{p,p}( C_{\rho}(t,x))}
\leq N\int_{B_{3\rho}}
 \int_{0}^{5\rho^{2}}\frac{1}{(|x|+\sqrt t)^{p}}\,dt \,dx
$$
$$
=N\rho^{d+2-p}\int_{|y|<3}
|y|^{2-p}\int_{0}^{5/|y|^{2}}\frac{dt}{
(1+\sqrt t)^{p }}dy. 
$$
The last factor of $N\rho^{d+2-p}$ is finite
if $p>2$ because
$$
\int_{0}^{\infty}\frac{dt}{
(1+\sqrt t)^{p }}<\infty 
$$
and $2-p>-d$. It is also 
obviously finite if $p\in(1,2]$.

In case $\max(|x|,\sqrt{|t|})\geq
2\rho$ for $(s,y)\in C_{\rho}(t,x)$, we have 
$$
|y|+\sqrt{|s|}\geq \max(|y|,\sqrt{|s|})
\geq \max(|x|,\sqrt{|t|}) 
-\max(|x-y|,\sqrt{|t-s|})\geq \rho.
$$
We see that in both cases
$$
\dashnorm g\|_{\L_{p,p}( C_{\rho}(t,x))}\leq
N\rho^{-1}.
$$

\end{remark}

\begin{remark}
                       \label{remark 5.14,5}
Theorem \ref{theorem 5.8,20} is also
true if we use the traditional way
to define the mixed-norm spaces integrating
first in $x$ and then in $t$. However, there are cases when \eqref{3.21.01} holds and
\begin{equation}
                              \label{5.14,4}
\sup_{\rho \leq1,C\in\bC_{\rho}}\rho 
\dashnorm b  \|_{ \L_{q\beta,p\beta}(C)} <\infty \quad
\text{and}\quad \sup_{\rho \leq q,C\in\bC_{\rho}}\rho 
\dashnorm b I_{C}\|_{ L_{p\beta,q\beta} } =\infty .
\end{equation}
Therefore, developing the theory of
solvability in $\L_{q,p,\beta}$ was worth
the effort. By the way, Minkowski's
inequality shows that \eqref{5.14,4}
is impossible if $q\geq p$.

To show an example of \eqref{5.14,4},
take $p,q,\beta>1$,  such that 
$$
\beta\leq \frac{d}{p}+\frac{2}{q},\quad 2\leq
q\beta<p\beta<d+1.
$$
  Set $p_{0}=\beta p,q_{0}=\beta q$ and
leq
  $b$ be such that
$$
|b(t,x)|=f(t,x):=I_{ t>0 }|x|^{-1}\Big| \frac{\sqrt t}{|x|}-1 \Big|^{-1/p_{0}}.
$$
Observe that for any constant $\lambda>0$
we have $f(\lambda^{2}t,\lambda x)=\lambda^{-1}f(t,x)$, implying that for any $C_{r}(t,x)$
\begin{equation}
                          \label{8.27.1}
\dashnorm f\|_{\L_{q_{0},p_{0}}(C_{r}(t,x))}=r^{-1}\dashnorm f\|_{\L_{q_{0},p_{0}}(C_{1}(t/r^{2},x/r))}.
\end{equation}

Next, 
$$
\int_{0}^{1}|f(t,x)|^{q_{0}}\,dt=|x|^{2-q_{0}}
\int_{0}^{1/|x|^{2}}|\sqrt{s}-1|^{-q/p}\,ds=:
|x|^{2-q_{0}}J(1/|x|^{2}).
$$
Here $q/p<1$, so that the singularity in $J$
is integrable and $J(t)\sim t^{1-q/(2p)}$
as $t\to\infty$, so that
$J(t)\leq N t^{1-q/(2p)}$. It follows that 
$$
\Big(\int_{0}^{1}|f(t,x)|^{q_{0}}\,dt\Big)^{1/q_{0}}
\leq N|x|^{1/p_{0}-1},
$$
where $N$ is independent of $x$, and 
since $1-p_{0}  >-d$, the right-hand side
to the power $p_{0}$ is summable over $B_{1}$.
Together with \eqref{8.27.1} this yield
\begin{equation}
                                \label{1.23.1}
\dashnorm f\|_{\L_{q\beta ,p\beta }(C_{r})}\leq Nr^{-1}\quad r>0.
\end{equation}

If $t>0$ and $\rho(t,x)\leq 3r$, then
$$
\dashnorm f\|_{\L_{q \beta,p \beta}(C_{r}(t,x))}\leq N
\dashnorm f\|_{\L_{q\beta ,p\beta }(C_{4r} )}\leq Nr^{-1}.
$$
In case $t>0$ and $\rho(t,x)> 3r$
use \eqref{8.27.1} which shows that
\begin{equation}
                                \label{8.27.4}
\dashnorm f\|_{\L_{q \beta,p \beta}(C_{r}(t,x))}=r^{-1}\dashnorm f\|_{\L_{q \beta,p \beta}(C_{1}(\tilde t,\tilde x))}
=Nr^{-1}\|f\|_{\L_{q \beta,p \beta}(C_{1}(\tilde t,\tilde x))},
\end{equation}
where $\tilde t=t/r^{2},\tilde x=x/r$, so that
$\rho(\tilde t,\tilde x)>3$. This case we split
into two subcases: A) $\tilde t>1$, B) $|\tilde
x|>2$.

In case A)
$$
\int_{\tilde t}^{\tilde t+1}|f(t,x)|^{q_{0}}\,dt=|x|^{2-q_{0}}
\int_{\tilde t/|x|^{2}}^{(\tilde t+1)/|x|^{2}}|\sqrt{s}-1|^{-q/p}\,ds 
$$
$$
=I_{\tilde t>|x|^{2}}...+I_{\tilde t\leq|x|^{2}}...
=J_{1}+J_{2}\leq J_{1}+
J(2).
$$
Furthermore, $J_{1}(\tilde t)$ has negative derivative for $\tilde t>|x|^{2}$. Therefore,
$$
J_{1}\leq |x|^{2-q_{0}}\int_{1}^{2/|x|^{2}}
|\sqrt{s}-1|^{-q/p}\,ds .
$$
It follows as above that
$$
\Big(\int_{\tilde t}^{\tilde t+1}|f(t,x)|^{q_{0}}\,dt\Big)^{1/q_{0}}
\leq N|x|^{1/p_{0}-1}+N,
$$
where the integrals over unit balls
of the right-hand side to the power $p_{0}$
are dominated by a finite constant.
This provides the desired estimate
of the left-hand side of \eqref{8.27.4}
in case A).

In case B) one can use the same arguments
augmented by the fact that on $B_{1}(\tilde x)$
we have $|x|\geq 1$.
 
 Finally, notice that for $t\in(0,1)$
$$
\int_{|x|\leq 1}f^{p_{0}}(t,x)\,dx
=N\int_{0}^{1}
\rho^{d-1-p_{0}}\Big| \frac{\sqrt t}{\rho}-1 \Big|^{-1 }\,d\rho=\infty.
$$
This proves the second relation in \eqref{5.14,4}.

\end{remark}
\begin{remark}
With $f$ from Remark \ref{remark 5.14,5} and
$|b(t,x)|=c f(t,x) $ with sufficiently small $c>0$
the equation $\partial_{t}u+\Delta u+b^{i}D_{i}u-\lambda u
=f$ for sufficiently large $\lambda$ has a unique
solution in $u\in\E^{1,2}_{q,q,\beta}$ with $\beta<2$
if $f$ is, say bounded function with compact support.
This is a particular case of Theorem \ref{theorem 5.8,20}.
By embedding theorems for Morrey-Sobolev spaces $u$ is bounded and H\"older continuous.
However, if you use the theory of solvability
in Sobolev spaces you will only be able to conclude
that $u\in W^{1,2}_{q}$ which does not guarantee
that $u$ is even bounded if $q<(d+2)/2$.
This shows an advantage of considering Morrey-Sobolev spaces.
\end{remark}

\begin{remark}
Example 3.8 of \cite{Kr_25} shows that
the smallness assumption on $\bar b$
is unavoidable.
\end{remark}

\begin{remark}
If $a,b$ are independent of $t$, Theorem
\ref{theorem 5.8,20} yields a result
containing the {\em corresponding\/}  results
in  \cite{CLV_96}, \cite{FHS_17},  \cite{Li_03},
\cite{MT_13} applied to equations in the 
{\em whole\/} space on account of the presence of $b$
which is not necessarily in $L_{d,\loc}$
and can be such that $|b|=c/|x|$.  

\end{remark}

\mysection{the case of main coefficients depending  only on $t$}
                     \label{section 8.23.1}

Here we start on our way of proving Theorem
\ref{theorem 5.8,20} following a quite different path from\cite{Kr_25}.
Fix a $\delta\in(0,1]$ and let $a=a(t)$
be an $\bS_{\delta}$-valued function on $\bR$.

\begin{lemma}
                   \label{lemma 3.29.1}
Let $a \in \bS_{\delta}$. Then for any $d\times d $ 
symmetric matrix $u$
we have
\begin{equation}
                              \label{3.29.2}
\{a,u\}:=a^{ij}a^{kr}u_{ik}u_{jr}\geq \delta^{2}\sum_{i,j}u_{ij}^{2},
\end{equation}
\begin{equation}
                              \label{3.29.3}
 (1-\delta^{2})^{2} \{a,u\}\geq \{a-\delta(\delta^{ij}),u\} \geq 0. 
\end{equation}
\end{lemma}

Proof. If $\lambda_{p},\ell_{p}$ are the eigenvalues and eigenvectors of $a$, then
the left-hand side of \eqref{3.29.2} is written
as
$$
\sum_{p,q}\lambda_{p}\lambda_{q}\Big(\sum_{i,j}
u_{ij}\ell_{p}^{i}\ell_{q}^{j}\Big)^{2}
\geq\delta^{2}\sum_{p,q} \Big(\sum_{i,j}
u_{ij}\ell_{p}^{i}\ell_{q}^{j}\Big)^{2}
=\delta^{2}u_{ij}\ell_{p}^{i}\ell_{q}^{j}
u_{kr}\ell_{p}^{k}\ell_{q}^{r}
$$
and the latter equals the right-hand side of
\eqref{3.29.2} because $\ell^{i}_{p}\ell^{k}_{p}=\delta^{ik},\ell^{j}_{q}\ell^{r}_{q}=\delta
^{jr}$.  

Since $a\in\bS_{\delta}$, $(1-\delta^{2})\lambda_{p}\geq \lambda_{p}-\delta\geq0$.
This yields \eqref{3.29.3}. \qed

Let $u\in C^{\infty}_{0}$. 
Set 
$$
-f=\partial_{t}u+a^{ij}(t)D_{ij}u.
$$
 Then
$$
\int_{\bR^{d+1}}f^{2}\,dz=I_{1}+ I_{2}+I_{3},
$$
where integrating by parts 
and observing that
$$
\int_{\bR^{d+1}}\partial_{t}u \delta^{ij}D_{ij}u\,dz=-(1/2)\int_{\bR^{d+1}}\partial_{t}(|Du|^{2})\,dz=0,
$$
we find
$$
I_{1}=\int_{\bR^{d+1}}|\partial_{t}u|^{2}\,dz,
$$
$$
I_{3}=\int_{\bR^{d+1}}a^{ij}D_{ij}ua^{kr}D_{kr}u\,dz=\int_{\bR^{d+1}}\{a,D^{2} u\}\,dz.
$$
$$
I_{2}=2\int_{\bR^{d+1}}\partial_{t}u a^{ij}D_{ij}u\,dz=2\int_{\bR^{d+1}}\partial_{t}u [a^{ij}-\delta\delta^{ij}]D_{ij}u\,dz 
$$
$$
\geq-2I_{1}^{1/2}\Big(\int_{\bR^{d+1}}\{a
-\delta (\delta^{ij}),D^{2} u\}\,dz\Big)^{1/2}
\geq-2(1-\delta^{2})I_{1}^{1/2}I_{3}^{1/2}.
$$
It follows that
$$
I_{1}+I_{3}\leq\delta^{-2}\int_{\bR^{d+1}}f^{2}\,dz, 
$$
\begin{equation}
                          \label{3.29,3}
\int_{\bR^{d+1}}\big(|\partial_{t}u|^{2}+\delta^{2}|D^{2}u|^{2}\big)\,dz\leq
\delta^{-2}\int_{\bR^{d+1}}f^{2}\,dz. 
\end{equation}

Introduce $B^{0,\infty}$ as the space
of function $u(t,x)$ on $\bR^{d+1}$
which are infinitely differentiable in $x$
for every $t$ with each derivative
locally bounded on $\bR^{d+1}$. By $B^{1,\infty}$ we mean the 
subspace of $B^{0,\infty}$ consisting of functions
$u(t,x)$ such that $\partial_{t}u\in B^{0,\infty}$. By $B^{i,\infty}_{0}$
we mean the subspaces of $B^{i,\infty} $  consisting of functions
with compact support, $i=0,1$.

 Note that 
  estimate \eqref{3.29,3} is also true  
for
  $u\in B^{1,\infty}$  such that, for some $\alpha>d/2$ and $\rho(t,x)=\sqrt{|t|}+|x|,\xi 
=\rho^{ \alpha-1},\eta=\rho^{\alpha}$ the function 
\begin{equation}    
                             \label{3.30.2}
\xi|u|+\eta|Du| 
\end{equation}
 is bounded. This is easily proved by
taking $u(t,x)\zeta(n^{2}t,nx)$ in place of $u$ with
$\zeta\in C^{\infty}_{0} $
such that $\zeta(0)=1$ and sending $n\to\infty$.

The following is just a reminder of standard
computations. 
Define 
$$
A_{t,t+ s}=\int_{t}^{t+s }a(r)\,dr
=s\int_{0}^{1}a(t+rs)\,dr,
\quad
\sigma_{t,t+ s}= s ^{-1/2}A^{1/2}_{t,t+s} ,
$$
recall that for $c_{d}=(4\pi)^{-d/2}$ 
$$
\sfp(t,x)=c_{d}t^{-d/2}e^{-|x|^{2}/(4t)}
I_{t>0}
$$
is the fundamental solution of the heat equation
and for $\lambda\geq0$ introduce
$$
R_{\lambda}f(t,x)=\int_{0}^{\infty}
\int_{\bR^{d}}e^{-\lambda s}f(t+s, x+
\sigma_{t,t+ s}y)\sfp(s,y)\,dyds
$$
$$
=e^{\lambda t}\int_{t}^{\infty}
\int_{\bR^{d}}e^{-\lambda s}f( s, x+
\sigma_{t,  s}y)\sfp(s-t,y)\,dyds 
$$
$$
=\int_{t}^{\infty}e^{-\lambda s}
\int_{\bR^{d}}f( s,y )\sfp(t, s,y-x)\,dyds,
$$
where 
$\sfp(t,s, x)=\sfp(s-t ,\sigma_{t, s}^{-1}x)\det \sigma_{t, s }^{-1}$.

\begin{remark}
                   \label{remark 3.30.2}
Observe that $\sigma_{t, s}\geq\delta^{1/2}$, which  implies that
for $n\geq0$ and $s>t$
\begin{equation}
                            \label{3.31.1}
|D^{n}\sfp(t,s, x)|\leq N(d,n,\delta)
\frac{1}{(s-t)^{(d+n)/2}}e^{-|x|^{2}\delta/(8  (s-t))}.
\end{equation}
Furthermore, from
$\sfp(t,s, x)=\sfp(s-t ,\sigma_{t, s}^{-1}x)\det \sigma_{t, s }^{-1}$ we get that the Fourier transform $\tilde p$ of $p$ satisfies
$$
\tilde p(t,s,\xi)=\exp(-(A_{t,s}\xi,\xi)),\quad
\partial_{t}\tilde p(t,s,\xi)
-(a(t)\xi,\xi)\tilde p(t,s,\xi)=0.
$$
It follows that
$$
\partial_{t} p(t,s,x)+a^{ij}(t)D_{ij}
p(t,s,x)=0 
$$
and \eqref{3.31.1} implies that for $n\geq0$ and $s>t$
\begin{equation}
                            \label{3.31,2}
|\partial_{t}D^{n}\sfp(t,s, x)|\leq N(d,n,\delta)
\frac{1}{(s-t)^{1+(d+n)/2}}e^{-|x|^{2}\delta/(8  (s-t))}.
\end{equation}

\end{remark} 

For $k,s,r>0,\alpha\in \bR $, and appropriate $f(t,x)$'s
on $\bR^{d+1}$
 define
$$
p_{\alpha,k}(s,r)=\frac{1}{s^{(d+2-\alpha)/2}}e^{-r^{2}/(ks)}I_{s>0}, 
$$
$$
P_{\alpha,k}f(t,x)=\int_{\bR^{d+1} }p_{\alpha,k}(s,|y|)f(t+s,x+y)\,dyds 
$$
$$
=\int_{t}^{\infty}\int_{\bR^{d} }p_{\alpha,k}(s-t,|y-x|)
f(s,y)\,dsdy.
$$  

\begin{remark}
                  \label{remark 5.31,1}
For $\alpha\leq d+2$ there is a constant $N=N(d,\alpha)$
such that
\begin{equation}
                        \label{5.31,1}
p_{\alpha,k}(s,|x|)\leq N\frac{1}{\rho
^{d+2-\alpha}(s,x)},
\end{equation}
where $\rho(s,x)=|x|+\sqrt{|t|}$.

Indeed, if $|x|\leq \sqrt{|s|}$, then
$$
p_{\alpha,k}(s,|x|)\leq 
\frac{1}{s
^{(d+2-\alpha)/2}}=\frac{1}{\rho
^{d+2-\alpha}(s,x)}\frac{\rho
^{d+2-\alpha}(s,x)}{s
^{(d+2-\alpha)/2}},
$$
where the last fraction is dominated
by $2^{d+2-\alpha}$. In case
$|x|\geq \sqrt{|s|}$,
$$
p_{\alpha,k}(s,|x|)\leq 
\frac{1}{|x|
^{d+2-\alpha}}\phi(|x|/\sqrt{|s|}),
$$
where $\phi(t)=t^{d+2-\alpha }e^{-t^{2}/k}$ is a bounded function on
$(0,\infty)$. In that case we get 
\eqref{5.31,1} again.

As a consequence of this estimate we
obtain that, if $f$ is a bounded function
with compact support, then 
\begin{equation}
                        \label{5.31,2}
|P_{\alpha,k}f(t,x)|
\leq N(1+\rho(t,x))^{-(d+2-\alpha)},
\end{equation}
where $N$ is independent of $(t,x)$.

\end{remark}

Introduce
$$
\cL^{0} u(t,x)=\partial_{t}u(t,x)+a^{ij}(t)D_{ij}u(t,x).
$$

The following two results are well known
and proved, for instance, by using the Fourier transform
with respect to $x$. Assertion (ii) in Theorem
\ref{theorem 3.30.1} follows from the fact that
$$
|R_{\lambda}f|\leq N(1+\rho )^{-d},\quad
|DR_{\lambda}f|\leq N(1+\rho )^{-(d+1)},
$$
which is a consequence of \eqref{5.31,2}.

\begin{theorem}
                     \label{theorem 3.30.1}
 Let $ f\in B^{0,\infty}_{0}$. Then

(i) we have $u:=R_{\lambda}f\in B^{1,\infty}$
and   
\begin{equation}
                             \label{3.30,3}
\cL^{0} u-\lambda u =-f   ;
\end{equation}
 
(ii) there is $\alpha>d/2$ such that the function \eqref{3.30.2} is bounded.
\end{theorem}

\begin{lemma}
                        \label{lemma 3.31.1}
Let $u\in C^{\infty}_{0}$, $\lambda
\geq0$. Then
$$
u=R_{\lambda}(\lambda u-\cL^{0} u).
$$
\end{lemma}
 
Define
$$
\osc_{C}g:=\dashint_{C}\dashint_{C}
|g(z_{1})-g(z_{2})|\,dz_{1}dz_{2},
\quad g^{\sharp}(t,x)=\sup_{\substack{C\in\bC,\\C\ni(t,x)}}\osc_{C}g .
$$
Also for $\beta\geq0$ define
$$
\bM_{\beta}g(t,x)=\sup_{\rho>0}\rho^{\beta}\sup_{\substack{C\in\bC_{\rho},\\C\ni(t,x)}}\dashint_{C}|g|\,dz,\quad \bM g=\bM_{0}g.
$$

Here is a particular case of
Lemma 2.2 of \cite{Kr_23}.

\begin{lemma}
                   \label{lemma 19.30.1}
For any $ \alpha<\beta $, $\beta\geq0$, and $k>0$, there exists a   constants $N$
  such that for any $f\geq0$, $\rho>0$,
 we have   on $C_{\rho}$ that
\begin{equation}
                     \label{1.17.20}
P_{\alpha,k}(I_{C^{c}_{2\rho}}f)  
\leq N \rho^{\alpha-\beta}\bM_{\beta}
f(0)  .
\end{equation}
\end{lemma}
 
We need this for the following.

\begin{theorem}
                        \label{theorem 3.31.1}
Let $f\in B^{0,\infty}_{ 0}$ and $u=R_{0}f$. Then
for any $\kappa\geq 2$ and $\rho>0$ 
$$
\osc_{C_{\rho}}D^{2}u\leq N(d,\delta)\kappa^{(d+2)/2}
\dashnorm f\|_{L_{2}(C_{\kappa\rho})}+
N(d,\delta)\kappa^{-1} \bM  f  (0).
$$
In particular, if $u\in C^{\infty}_{0}$ then 
$$
\osc_{C_{\rho}}D^{2}u\leq N(d,\delta)\kappa^{(d+2)/2}
\dashnorm \cL^{0} u\|_{L_{2}(C_{\kappa\rho})}+
N(d,\delta)\kappa^{-1} \bM (\cL^{0} u)  (0). 
$$

\end{theorem}

Proof. Observe that the case of arbitrary $\rho>0$ is reduced to $\rho=1$ by using
scale changes. In  case $\rho=1$ 
take $\zeta \in C^{\infty}_{0}$
such that $\zeta=1$ in $C_{\kappa}$, $\zeta=0$ 
outside $C_{ 2\kappa}^{c}\cap \bR^{d+1}_{0}$, $0\leq\zeta\leq1$ and  
note that for   $g=f\zeta,h=f(1-\zeta)$
 and $(G,H)=R_{0}(g,h)$   we have
$u=G+H$, and for $m=0,1,n\geq0$    
 $$
 \partial_{t} ^{m}D^{n}
H(t,x)= \int_{C_{2}^{c}}
 \partial_{t} ^{m}D^{n}
\sfp(t, s,y-x)h( s, y)\,dsdy.   
$$

By using   \eqref{3.31.1} we get that in $C_{1}$  
$$
|\partial_{t}D^{2}H |+|D^{3}H |
\leq P_{-2,\delta/8 }|h|+P_{-1,\delta/8 }|h|, 
$$
which by Lemma \ref{lemma 19.30.1} 
(with $\beta=0$) implies that 
\begin{equation}
                         \label{1.26,1}
|\partial_{t}D^{2}H(t,x)|+|D^{3}H(t,x)|
\leq N\kappa^{-1}\bM h(0).
\end{equation}
It follows  that 
$$
\int_{C_{1}}\int_{C_{1}}
|D^{2}H(z_{1})-D^{2}H(z_{2}|\,dz_{1}dz_{2}\leq N
\kappa^{-1}\bM f(0). 
$$

Regarding $G$ we have in light of   \eqref{3.29,3} 
 and Theorem \ref{theorem 3.30.1}
that 
$$
\int_{C_{1}}\int_{C_{1}}
|D^{2}G(z_{1})-D^{2}G(z_{2})|\,dz_{1}dz_{2}\leq N\int_{C_{1}} 
|D^{2}G(t,x) |\,dxdt 
$$
$$
\leq N\|g\|_{L_{2}(\bR^{d+1})}\leq N
\kappa^{(d+2)/2}\dashnorm f\|_{L_{2}(C_{\kappa })}.
$$
This proves the first assertion of the theorem. 
The second one is a consequence of the first one and Lemma \ref{lemma 3.31.1}.  \qed

Let $p,q\in(1,\infty)$.
 Recall the definition of
$\L_{q,p}$  and let
the space $\W^{1,2}_{q,p}(Q)$ be defined as 
the collection of functions $u(t,x)$ on $Q$
such that $u,Du,D^{2}u,\partial_{t} u
\in \L_{q,p}(Q)$. The space $\W^{1,2}_{q,p}(Q)$ is provided with a natural norm. We set
$\W^{1,2}_{q,p} =\W^{1,2}_{q,p}(\bR^{d+1})$.
\begin{theorem}
                         \label{theorem 3.31.2}
There is a constant $N=N(d,\delta,p,q)$
such that for any $u\in \W^{1,2}_{q,p}$
\begin{equation}
                               \label{4.1,5}
\|\partial_{t}u,D^{2}u\|_{\L_{q,p}}
\leq N\|\cL^{0} u\|_{\L_{q,p}}.
\end{equation}
\end{theorem}

Proof. By Lemma \ref{lemma 3.31.1} and Theorem \ref{theorem 3.31.1} with $\kappa=2$ for $u\in C^{\infty}_{0}
(\bR^{d+1})$ we have
$$
(D^{2}u)^{\sharp}\leq N(d,\delta)(\bM((\cL^{0} u)^{2}))^{1/2}, 
$$
which by the Muckenhoupt weights versions of the Fefferman-Stein and Hardy-Littlewood theorems,
proved by Hongjie Dong and Doyoon Kim in \cite{DK_18}, yields
for $p>2 $ and $w\in A_{p}$ ($A_{p}$ is the
Muckenhoupt class of weights) 
$$
\int_{\bR^{d+1}}|D^{2}u|^{p}w\,dxdt
\leq N\int_{\bR^{d+1}}|\cL^{0} u|^{p}w\,dxdt. 
$$
After that the estimate for $D^{2}u$ 
 follows from
another Dong-Kim theorem found in \cite{DK_18}. 
The term $\partial_{t}u$ as usual
is estimated from $\partial_{t}u=\cL u-a^{ij}D_{ij}u$. \qed

Estimate \eqref{4.1,5} for $p=q$ (and variable
$a$) was first obtained in \cite{Jo_64},
for $p\ne q$ it is proved in \cite{Kr_01}.
The ideas used here are quite different
from both \cite{Jo_64} and \cite{Kr_01}
and are borrowed from \cite{DK_18},
where mixed norm estimate are obtained
with the integration on $x$ and $t$ in any order.

Next, we introduce the term $\lambda u$
by using the well-known Agmon's method 
(see, for instance, the proof of Lemma 6.3.8
of \cite{Kr_08}) and adding 
an interpolation inequality, allowing
to estimate the norms of $Du$ through
the norms of $u$ and $D^{2}u$ (see Remark \ref{remark 8.22.3} below), conclude that
the a priori estimate \eqref{4.1,1}
with $s=-\infty$
is true for $u\in C^{\infty}_{0}$.
For any $u\in \W^{1,2}_{q,p}$ it holds
because $C^{\infty}_{0}$ is dense
in $\W^{1,2}_{q,p}$.

Define
$$
\bR^{d+1}_{s}=(s,\infty)\times\bR^{d}.
$$

\begin{theorem}
                         \label{theorem 4.2,1}
There is a constant $N=N(d,\delta,p,q)$
such that for any $s\in[-\infty,\infty)$,
 $u\in \W^{1,2}_{q,p}(\bR^{d+1}_{s})$ and
$\lambda\geq0$ 
\begin{equation}
                               \label{4.1,1}
\|\lambda u,\sqrt\lambda Du,\partial_{t}u,D^{2}u\|_{\L_{q,p}(\bR^{d+1}_{s})}
\leq N\|\cL^{0} u-\lambda u\|_{\L_{q,p}(\bR^{d+1}_{s})}.  
\end{equation}
Furthermore, for any $\lambda>0$ and
$f\in \L_{q,p}(\bR^{d+1}_{s})$ there is a unique
$u\in \W^{1,2}_{q,p}(\bR^{d+1}_{s})$ such that
$\cL^{0} u-\lambda u=f$ in $\bR^{d+1}_{s}$.  
\end{theorem}

The existence result for $s=-\infty$
follows from \eqref{4.1,1} and Theorem \ref{theorem 3.30.1}
in which, as is easy to see, $u\in\W^{1,2}_{q,p}(\bR^{d+1})$. For $s>-\infty$ estimate
\eqref{4.1,1} and the existence part follow from the case $s=\infty$ and uniqueness
when $s>-\infty$. In turn
the uniqueness of $\W^{1,2}_{q,p}(\bR^{d+1}_{s})$-solutions (the causality property) follows from the representation of solutions in Lemma 
\ref{lemma 3.31.1}, which   extends
to $u\in \W^{1,2}_{q,p}(\bR^{d+1}_{s})$
by continuity and \eqref{4.1,1}.

\begin{remark}
                    \label{remark 8.22.3}
Above and a few times below we use interpolation inequalities of the following type.
Take $a=(\delta^{ij})$ and notice that  according to
\eqref{3.31.1} 
$$
|DR_{\lambda}f(t,x)|\leq N
\int_{0}^{\infty}\int_{\bR^{d}}e^{-\lambda s}
\frac{|f(s+t,x+y)|}{s^{(d+1)/2}}e^{-|y|^{2} /(8 s)}\,dyds,
$$
which estimates  $|DR_{\lambda}f|$ in term of the ``sum'' with respect to $(s,y)$
with weights depending only on $(s,y)$
of functions $|f(s+t,x+y)|$ as functions of $(t,x)$. Minkowski's inequality that the norm
of a sum is less than the sum of norms
now leads to the first assertion in the following. The second assertion follows
from Lemma \ref{lemma 3.31.1}.

\begin{lemma}
                      \label{lemma 3.31.4}
Let $s\in[-\infty,\infty)$ and let $\frB$ be a Banach space of
functions on $\bR^{d+1}_{s}$ such that if $f\in \frB$
and $(s,y)\in\bR_{0}^{d+1}$, then $f_{s,y}\in\frB$ and
$\|f_{s,y}\|_{\frB}\leq \|f\|_{\frB}$, where
 $f_{s,y}(t,x)=f(s+t,x+y)$. Also assume that,  
if $g\in\frB$ and a measurable $f$ satisfies $|f|\leq g$, then $f\in\frB$ and  $\|f\|_{\frB}\leq \|g\|_{\frB}$. Then for any $\lambda>0$
and $f\in\frB$ we have 
$$
 \lambda^{1/2}\|DR_{\lambda}f\|_{\frB}\leq N(d)\|
f\|_{\frB}.
$$
Furthermore, if $C^{\infty}_{0}
\subset \frB$ then for any
$u\in C^{\infty}_{0}$ and $\lambda>0$ 
$$
\|Du\|_{\frB}\leq N\lambda^{-1/2}
\|\partial_{t}u+\Delta u-\lambda u\|_{\frB}
\leq N\big(\lambda^{-1/2}
\|\partial_{t}u,D^{2} u \|_{\frB}
+\lambda^{1/2}
\| u\|_{\frB}\big), 
$$
where $N=N(d)$.

\end{lemma}
\end{remark}

\mysection{Variable main coefficients and Sobolev spaces}
 \label{section 4.29.1}
 
Now we are going to extend Theorem
\ref{theorem 3.31.2} to operators with 
$\bS_{\delta}$-valued $a=a(t,x)$.
Define 
$$
a^{\shharp} =a^{\shharp}_{\infty}=
 \sup
_{ C\in\bC }\dashint_{C}|a(t,x)-\tilde a_{C}(t)|\,dxdt,\quad \tilde a_{C}(t)=\dashint_{C}a(t,x)\,dxds.
$$

Observe that if $a$
is independent of $x$, then $a^{\shharp}  =0$. 

Set
$$
\cL_{0}  u=\partial_{t}u+a^{ij}D_{ij}u.
$$

For  $ C\in\bC$  and
$$
 \cL_{ C} u:= \partial_{t}u+
\tilde a^{ij}_{ C}D_{ij}u,\quad f=\cL_{0}u
$$
  we have 
$$
f= \cL_{ C} u
+[(a^{ij}- \tilde a^{ij}_{C})D_{ij}u],
$$
 whereas
for $q>2$
$$
\dashnorm (a^{ij}- \tilde a^{ij}_{ C})D_{ij}u
\|_{L_{2} C )}\leq
N(d,\delta, q)a^{\shharp} 
\big(\bM(|D^{2} u|^{q})\big)^{1/q}. 
$$
This and Theorem \ref{theorem 3.31.1} lead  to the following.

\begin{lemma}  
                  \label{lemma 4.2,1}

For any $q> 2$,
$\kappa\geq 2$, and $u\in C^{\infty}_{0}$   we have 
$$
( D^{2}u)^{\sharp}
\leq N(d,\delta,q,\kappa)\Big(\bM(|\cL_{0} u|^{q})\big)^{1/q} 
$$
\begin{equation}
                         \label{4.2,1}
+\big(N(d,\delta,q,\kappa)a^{\shharp} +N(d,\delta )\kappa^{-1}\big)\big(\bM(|D^{2} u|^{q})+\bM(|\partial_{t}u|^{q})\big)^{1/q}\Big).
\end{equation}
\end{lemma}  
 
Since $\partial_{t}u=\cL_{0}u-a^{ij}D_{ij}u$,
similar estimate holds for $\partial_{t}u$.
Therefore, the Fefferman-Stein and Hardy-Littlewood allow us to do the next step
for $p>2$ as in the case of Theorem \ref{theorem 3.31.2}.
\begin{lemma}   
                        \label{lemma 4.2,6}
For any $p>1$, $w\in A_{p}$, $\kappa\geq2 $,
and $u\in C^{\infty}_{0}$
$$
\|D^{2}u\|_{L_{p}(w)}^{p}+\|\partial_{t}u\|_{L_{p}(w)}^{p}
\leq N([w]_{A_p})\Big(N(d,\delta,p,\kappa)
\|\cL_{0} u\|_{L_{p}(w)}^{p} 
$$
\begin{equation}
                      \label{4.2,7}
+\big(N(d,\delta,p,\kappa)a^{\shharp} +N(d,\delta )\kappa^{-1}\big)^{p}\Big( \|D^{2}u\|_{L_{p}(w)}^{p}+\|\partial_{t}u\|_{L_{p}(w)}^{p}\big)\Big). 
\end{equation}
\end{lemma}

We will use a particular $A_{1}$-weight.
By the way, recall that $A_{1}$-weights
are also $A_{p}$-weights if $p>1$.

\begin{lemma}
                          \label{lemma 1.29.3}
Define $w_{\alpha}=(|x|+ \sqrt{|t|})^{-\alpha}$. Then for any $\alpha\in[0,d+2)$ the function $w_{\alpha}$
 is an $A_{1}$-weight, that is, there
is a constant $N$ such that $\bM w_{\alpha}
\leq Nw_{\alpha}$.
\end{lemma}

This lemma is proved by routine estimates.
The most important property of $w_{\alpha}$
for us
is the following.

\begin{lemma}
                        \label{lemma 1.30.2}
Let $\alpha\geq 0,\alpha+\beta>d+2 $. Then
for any $f\geq0$  
\begin{equation}
                            \label{1.30.4}
\int_{\bR^{d+1}}f(w_{\alpha}\wedge 1)\,dyds\leq  
\bM_{\beta}f(0).
\end{equation}

\end{lemma}

Proof. Without losing generality we suppose that $f$ is bounded and has compact support. Set 
$$
Q_{1}=\{(s,y):|y|\geq\sqrt {|s| }\},\quad Q_{2}=\{(s,y):|y|\leq \sqrt{ |s|}\}.
$$
 Note that
on $Q_{1}$ we have $w_{\alpha}(s,y)\leq N |y|^{-\alpha}$
and
$$
\int_{Q_{1}}f(w_{\alpha}\wedge 1)\,dyds\leq N\int_{\bR^{d}}
(1\wedge|y|^{-\alpha})\int_{-|y|^{2}}^{|y|^{2}}
f(s,y)\,dsdy 
$$
$$
\leq N\int_{0}^{\infty}(1\wedge r^{-\alpha})
\Big(\frac{\partial}{\partial r}\int_{B_{r}}
\int_{- |y|^{2}} ^{ |y|^{2}} 
f(s,y)\,ds\,dy\Big)\,dr
$$
$$
=N\int_{1}^{\infty}r^{-\alpha-1}
 \int_{B_{r}}
\int_{- |y|^{2}} ^{ |y|^{2}} 
f(s,y)\,ds\,dy \,dr
$$
$$
\leq N\bM_{\beta}f(0)\int_{1}^{\infty}r^{-\alpha-\beta+(d+2)-1}\,dr\leq N\bM_{\beta}f(0).
$$

Next, on $Q_{2}$ we have $w_{\alpha}(s,y)\leq N |s|^{-\alpha/2}$
and 
$$
\int_{Q_{2}}f(w_{\alpha}\wedge 1)I_{s>0}\,dyds\leq N\int_{0}^{\infty}
(1\wedge s^{-\alpha/2})\int_{B_{\sqrt s}}
f(s,y)\,dyds
$$
$$
 \leq N\int_{0}^{\infty}
(1\wedge s^{-\alpha/2})\Big(\frac{\partial}{\partial s}\int_{0}^{s}\int_{B_{\sqrt s}}
f(t,y)\,dtdy\Big)ds
$$
$$
\leq N\int_{1}^{\infty}
  s^{-\alpha/2-1}  \int_{0}^{s}\int_{B_{\sqrt s }}
f(t,y)\,dtdy ds
$$
$$
\leq N\bM_{\beta}f(0)\int_{1}^{\infty}
  s^{-\alpha/2+(d+2)/2-\beta/2-1}\,ds
\leq N\bM_{\beta}f(0).
$$
Similarly we estimate the part of the integral
of $f(w_{\alpha}\wedge1)$ over $Q_{2}\cap\{s<0\}$.   \qed

Take $p,q\in(1,\infty)$ and
   $\beta> 0$ and introduce the
Morrey space $\dot \E_{q,p,\beta} $
as the set of $g\in  \L_{q,p,\loc}$ 
such that  (note $\rho>0$)
\begin{equation}
                             \label{2.16.50}
\|g\|_{\dot \E_{q,p,\beta} }:=
\sup_{\rho >0,C\in\bC_{\rho}}\rho^{\beta}
\dashnorm g  \|_{ \L_{q,p}(C)} <\infty .
\end{equation}  
Define
$$
\dot \E^{1,2}_{q,p,\beta} =\{u:u,Du,D^{2}u,
\partial_{t}u\in \dot \E_{q,p,\beta} \},
$$
where $Du,D^{2}u,
\partial_{t}u$ are the Sobolev
 derivatives,
and 
provide $\dot \E^{1,2}_{q,p,\beta} $ with an obvious norm.
If $p=q$ we drop $q$ in the {\em notation\/} of these spaces.

Now we can make the first crucial step
towards proving Theorem \ref{theorem 5.8,20}.

\begin{theorem}
                \label{theorem 4.30,1}
Let $p \in(1,\infty)$, $ \beta>0$. Then
there is a constant $\hat a=\hat a(d,\delta,p,\beta) >0$ and a constant
$N=N(d,\delta,p,\beta)$ such that if
\begin{equation}
                       \label{4.30,2}
a^{\shharp}\leq \hat a,
\end{equation}
then for any $\lambda\geq0$ and
$u\in \dot E^{1,2}_{p,\beta}$
\begin{equation}
                       \label{4.30,3}
\|\lambda u,\sqrt\lambda Du,\partial_{t}u,D^{2}u\|_{\dot E_{p,\beta}
 }\leq N\|\cL_{0} u-\lambda u\|
_{\dot E_{p,\beta} }.
\end{equation}
Furthermore, for any $\lambda>0$
and $f\in \dot E_{p,\beta} $
there exists a unique solution 
$u\in \dot E^{1,2}_{p,\beta}$
of the equation $\cL_{0} u-\lambda u=f$.
\end{theorem}

Proof. Just in case, note that, if $\beta>d/p$,
our assertions are trivial because then $\dot E_{p,\beta}$ consists of only zero function. Then to prove \eqref{4.30,3}
observe that the case of arbitrary $\lambda\geq0$
is reduced to $\lambda=0$ by Agmon's method and interpolation. After that it is not hard to see
that we may assume that $u$ has
compact support, so that $u\in W^{1,2}_{p}$.

Now take $\alpha\in[0,d+2)$ such that $\alpha+p\beta>d+2$
and let $w =w_{\alpha}\wedge 1$,
so that $w$ is an $A_{1}$-weight owing to Lemma \ref{lemma 1.29.3}.
Approximating $u$ by smooth functions
in $W^{1,2}_{p}$ we conclude that \eqref{4.2,7}
holds for our $u$.

 Then we borrow an idea from \cite{CF_88}. By using
Lemma \ref{lemma 1.30.2}, we get from
\eqref{4.2,7} and the Muckenhoupt theorem 
and H\"older's inequality
that
with $w=w_{\alpha}\wedge 1$ we have
$$
\Big(\dashint_{C_{1}}|D^{2}u|^{p}\,dxdt\Big)^{1/p}\leq N
\|D^{2}u\|_{L_{p}( w)}
 \leq N\|\cL_{0} u \|
_{\dot E_{p,\beta}(\bR^{d+1})} 
$$
$$
+\big(N(d,\delta,p,\kappa,\beta)a^{\shharp} +N(d,\delta,\beta )\kappa^{-1}\big)\|\partial_{t}u,D^{2}u\|_{\dot E_{p,\beta}
(\bR^{d+1})}. 
$$

Shifting $w$, using the self-similarity and also
treating $\partial_{t}u$ in the same way give      
$$
\|\partial_{t}u,D^{2}u\|_{\dot E_{p,\beta}
(\bR^{d+1})}\leq N\|\cL_{0} u \|
_{\dot E_{p,\beta}(\bR^{d+1})} 
$$
$$
+\big(N(d,\delta,p,\kappa,\beta)a^{\shharp} +N(d,\delta,\beta )\kappa^{-1}\big)\|\partial_{t}u,D^{2}u\|_{\dot E_{p,\beta}
(\bR^{d+1})}, 
$$
which easily yields \eqref{4.30,3}
with $\lambda=0$.

The solvability part of the theorem,
as usual, is proved by the method of continuity starting from the heat equation when
the explicit formulas for solutions are available.
\qed

Theorem \ref{theorem 4.30,1} is derived from 
Lemma \ref{lemma 4.2,6}. This lemma also allows us to get the following result by using the same arguments that lead to Theorem~\ref{theorem 3.31.2}.

\begin{lemma}
                    \label{lemma 4.2,2}
For any $p,q\in( 1,\infty)$,
$\kappa\geq 2$, and $u\in C^{\infty}_{0}$  we have 
$$
\|\partial_{t}u,D^{2}u\|_{\L_{q,p}}\leq
  N_{0}(d,\delta,p,q,\kappa) \|\cL_{0} u\|_{\L_{q,p}} 
$$
\begin{equation}
                         \label{4.2,2}
+\big(N_{1}(d,\delta,p,q,\kappa)a^{\shharp} +N_{2}(d,\delta, p,q)\kappa^{-1}\big)\|
\partial_{t}u,D^{2}u\|_{\L_{q,p}}.
\end{equation}
 
\end{lemma}

\begin{corollary}
                 \label{corollary 4.2,1}
Let $p,q\in( 1,\infty)$ and assume that
\begin{equation}
                         \label{4.2,3}
 a^{\shharp} \leq (1/4)
N_{1}^{-1}\big(d,\delta,p,q,4
N_{2} (d,\delta, p,q)\big).
\end{equation}
Then for any $u\in C^{\infty}_{0}$  we have
\begin{equation}
                                 \label{4.3,3}
\|\partial_{t}u,D^{2}u\|_{\L_{q,p}}\leq
  2N_{0}(d,\delta,p,q,4
N_{2} (d,\delta, p,q)) \|\cL_{0}u\|_{\L_{q,p}}, 
\end{equation}
where $N=N(d,\delta, p,q )$.
\end{corollary}

We restate this corollary as follows.

\begin{theorem}
                    \label{theorem 4.3,1}
Let $p,q\in(1,\infty) $. Then
there are   constants $\hat a =\hat a(d,\delta,p,q)>0$ and $N=N(d,\delta, p,q)$ such that if 
\begin{equation}
                         \label{4.3,5}
 a^{\shharp} \leq \hat a,
\end{equation}
then for any $u\in C^{\infty}_{0}$
we have
\begin{equation}
                         \label{4.2,5}
\|\partial_{t}u,D^{2}u\|_{\L_{q,p}}\leq
  N  \|\cL_{0} u\|_{\L_{q,p}}. 
\end{equation}

\end{theorem}

Next, as in the case of Theorem \ref{theorem 4.2,1} we introduce $\lambda\geq0$  to see that
for $\lambda\geq0$ and $u\in C^{\infty}_{0}$
$$
\|\partial_{t}u,D^{2}u,\sqrt\lambda Du,\lambda u\|_{\L_{q,p}}\leq
N\|\cL_{0} u-\lambda u\|_{\L_{q,p}}. 
$$
In this way we arrive at the first assertion in the following theorem
if  $u\in C^{\infty}_{0}$
and $s=-\infty$.   

\begin{theorem}
                 \label{theorem 4.3,3}
Let $p,q\in( 1,\infty)$, $s\in[-\infty,\infty)$ and assume \eqref{4.3,5}.  
Then there exists $N =N (d,\delta,p,q)$   such that for any $u\in \W^{1,2}_{q,p}
(\bR^{d+1}_{s})$
and $\lambda\geq 0$ 
\begin{equation}
                            \label{4.3,4}
\|\partial_{t}u,D^{2}u,\sqrt\lambda Du,\lambda u\|_{\L_{q,p}(\bR^{d+1}_{s})}\leq
N\|\cL_{0} u-\lambda u\|_{\L_{q,p}(\bR^{d+1}_{s})}. 
\end{equation}
Furthermore, for any $f\in \L_{q,p}(\bR^{d+1}_{s})$
and $\lambda>0$
there exists a unique $u\in \W^{1,2}_{q,p}(\bR^{d+1}_{s})$
satisfying in $\bR^{d+1}_{s}$
$$
\cL_{0} u-\lambda u+f=0.
$$
\end{theorem}

 For $s=-\infty$ the a priori estimate \eqref{4.3,4} extends
from $C^{\infty}_{0}$ to $\W^{1,2}_{q,p}$
due to the denseness of the former in the latter and the solvability result is obtained by the method of continuity starting from
Theorem  
\ref{theorem 4.2,1}. General $s$ is treated as in the proof of 
Theorem  
\ref{theorem 4.2,1} by observing that
the causality property is preserved under
the method of continuity. 

Here is a useful
interior estimate.

\begin{theorem}
                 \label{theorem 4.29,3}
Let $p,q\in( 1,\infty)$ and assume \eqref{4.3,5}.
Let $0<r<\rho<\infty$,    
$u\in \W^{1,2}_{q,p}(C_{\rho})$.
Set
$f:=\cL_{0} u $.
 Then   
 \begin{equation}
                                                 \label{06.5.25.20}
 \|\partial_{t}u, D^{2}u\|_{\L_{q,p}(C_{r})}
\leq N\big(
\|f\|_{\L_{q,p}(C_{\rho})} 
+(\rho-r)^{-2}\|u-c\|_{\L_{q,p}(C_{\rho})}\big), 
\end{equation} 
where $N=N(d,\delta,p,q)$ and $c=c(x)$ is any affine function.

\end{theorem}

Proof. We follow the proof of Lemma 2.4.4
of \cite{Kr_08}.   
For obvious reasons we may assume that $u\in C^{1,2}(\bar{C}_{\rho})$ and $c=0$.
Then, let $\chi(s)$ be an infinitely differentiable function
on $\bR$ such that $\chi(s)=1$ for $s\leq0$ and $\chi(s)=0$
for $s\geq1$. For $m=0,1,2,...$ introduce  ($r_{0}=r$)
 $$
r_{m}=r+(\rho-r)\sum_{j=1}^{m}2^{-j},\quad
\xi_{m}(x)=\chi\big(2^{m+1}(\rho-r)^{-1}(|x|-r_{m}) \big),
$$
 $$
\eta_{m}(t)=\chi\big(2^{2m+2}(\rho-r)^{-2}(t-r^{2}_{m})\big),
\quad \zeta_{m}(t,x)=\xi_{m}(x)\eta_{m}(t).    
$$ 
As is easy to check, for
$$
C(m)=C_{r_{m}}=(0,r_{m}^{2})\times B_{r_{m}}
 $$
it holds that
 $$
\zeta_{m}=1\quad\text{on}\quad C(m),\quad\zeta_{m}=0\quad\text{on}
\quad C_{\rho}\setminus C(m+1).
 $$
Also (observe that $N2^{m+1}=N_{1}2^{m}$ with $N_{1}=2N$) 
 $$
|D\zeta_{m }|\leq N2^{m}(\rho-r)^{-1},\quad
|\partial_{t}\zeta_{m}|\leq N2^{2m}(\rho-r)^{-2},
  $$
\begin{equation}
                                                                \label{06.5.25.1}
|D^{2}\zeta_{m}|\leq N2^{2m}(\rho-r)^{-2}.
\end{equation} 

Next, the function $\zeta_{m}u$ is in $\W^{1,2}_{q,p}(\bR^{d+1}_{0})$  
and satisfies 
 $$
\cL_{0}(\zeta_{m} u) =\zeta_{m} f
+u\cL_{0}\zeta_{m} +
2a^{ij}D_{i}\zeta_{m}D_{j}u .
 $$
By Theorem \ref{theorem 4.3,1}
and the above-mentioned properties of $\zeta_{m}$
$$
A_{m} :=\|\partial_{t}(\zeta_{m} u), D^{2}(\zeta_{m}u)\|_{\L_{q,p}(Q(m+1))}
             \label{06.5.25.2} 
 \leq N\|f\|_{\L_{q,p}(C_{\rho})} +N2^{2m}(\rho-r)^{-2}I+NJ_{m},
$$ 
where
 $$
I=\|u \|_{\L_{q,p}(C_{\rho})},
\quad J_{m}=\| 
2a^{ij}D_{i}\zeta_{m}D_{j}u\|_{\L_{q,p}(C_{\rho})}
 $$ 
 $$
 \leq  N2^{m}(\rho-r)^{-1}\|Du \|_{\L_{q,p}
(Q(m+1))}.
 $$
By interpolation inequalities  for  
any $\varepsilon>0$
$$
\|Du \|_{\L_{q,p}
(Q(m+1))}\leq \|D(\zeta_{m+1}u) \|_{\L_{q,p}
(\bR^{d+1}_{0})} 
$$
\begin{equation}
                       \label{8.12.3}
\leq \varepsilon (\rho-r) 2^{- m}
A_{m+1}+N\varepsilon^{-1}(\rho-r)^{-1}
2^{ m}I. 
\end{equation}

It follows that for any $\varepsilon\in(0,1]$
$$
A_{m}\leq N\|f\|_{\L_{q,p}(C_{\rho})} +N
\varepsilon^{-1}2^{2m}(\rho-r)^{-2}I
+\varepsilon A_{m+1}. 
$$
 
Now we take $\varepsilon=1/8$ and get 
 \begin{gather}\nonumber
\varepsilon^{m}A_{m} 
\leq N\varepsilon^{m}\|f\|_{\L_{q,p}(C_{\rho})}+N\varepsilon^{m-1}
  2^{2m}(\rho-r)^{-2} I+
\varepsilon^{m+1} A_{m+1},
\\ 
                                                                \label{06.5.25.3}
A_{0}+\sum_{m=1}^{\infty}\varepsilon^{m}A_{m} 
\leq N \|f\|_{\L_{q,p}(C_{\rho})}+N (\rho-r)^{-2} I +
\sum_{m=1}^{\infty}\varepsilon^{m}  A_{m }.
\end{gather} 
Here the series converges because owing to \eqref{06.5.25.1},
 $$
A_{m}\leq N(1+4^{m}(\rho-r)^{-2})\|u\|_{\W^{1,2}_{q,p}(C_{\rho})}.
 $$
Therefore upon cancelling like terms in \eqref{06.5.25.3},
we see that $A_{0}$ is less than the right-hand side of \eqref{06.5.25.20}.
Since its left-hand side is obviously less than $A_{0}$, the theorem is proved.
\qed

\begin{remark}
                   \label{remark 8.12.1}
Similarly we get that $A_{1}$
is dominated by the right-hand side of 
\eqref{06.5.25.20} and this along with
\eqref{8.12.3} shows that 
$\|  D  u\|_{\L_{q,p}(C_{r})}$
is dominated by the right-hand side of 
\eqref{06.5.25.20} as well.

\end{remark}

To proceed further take $p,q\in(1,\infty)$,
  $r\in(1,p\wedge q]$ and find an integer $n=n(d,p,q)\geq 2$ and $t_{0}=1>t_{ 1}>...> t_{n}=0$ such that
\begin{equation}
                            \label{8.10.1}
\frac{d}{p_{k+1}}+\frac{2}{q_{k+1}}\leq\frac{d}{p_{k }}+\frac{2}{q_{k }}+1,
\end{equation}
where
$$
p_{k}=r+t_{k}(p-r ),\quad
q_{k}=r+t_{k}(q-r).
$$
Note that $t_{k}$ and $p_{k},q_{k}$
decrease as $k$ increases, 
$(p_{0},q_{0})=(p,q)$, $(p_{n},q_{n})=(r,r)$.

The following result provides a key to extending
Theorem 7.27 of \cite{Kr_25_2} to $\beta\in(1,\infty)$
instead of $\beta\in(1,2)$.
\begin{theorem}
                        \label{theorem 8.11.1}
Suppose that condition \eqref{4.3,5}
is satisfied with 
\begin{equation}
                            \label{8.22.4}
\hat a=\hat a_{0}(d,\delta,p,q,r):=\min_{k=0,...,n}
\hat a(d,\delta,
p_{k},q_{k}).
\end{equation}
Let $u\in \W^{1,2}_{q,p}(C_{2})$
and $\cL_{0} u=0$ in $C_{2}$. Then
\begin{equation}
                     \label{8.11.2}
\|u\|_{\L_{q,p}(C_{1})}
\leq N\|u\|_{\L_{r}(C_{2})},
\end{equation}
 where $N$ depends only on $d,\delta,
p_{k},q_{k}$, $k=0,1,...,n$.

\end{theorem}

Proof. Let $\chi(s)$ be an infinitely differentiable function
on $\bR$ such that $\chi(s)=1$ for $s\leq0$ and $\chi(s)=0$
for $s\geq1/2$. For $k=0,1,...,n$ set
$$
\xi_{ k}(x)=\chi\big(n (|x|-1-k/n)\big) ,\quad
\eta_{ k}(t)=\chi\big(n^{2} (t-(1+k/n)^{2})\big) ,
$$
$$
\zeta_{ k}(t,x)=\eta_{ k}(t)\xi_{ k}(x),\quad \rho_{k}=1+k/n.
$$
Observe that $C_{\rho_{0}}=C_{1},C_{\rho_{n}}=C_{2}$
and
 $\zeta_{k}=1$ on $C_{\rho_{k}}$, $\zeta_{k}=0$ in $C_{2}
\setminus C_{\rho_{k+1/2}}$.

By Theorem 10.2 of \cite{BIN_75} for any $k
=0,1,...,n-1$
$$
\|u\|_{\L_{q_{k},p_{k}}(C_{\rho_{k}})}
\leq\|u\zeta_{ k}\|_{\L_{q_{k},p_{k}}(C_{\rho_{k+1}})}
$$
\begin{equation}
                      \label{8.12.1}
\leq N\big(\|\partial_{t}
(u\zeta_{ k}),D^{2}(u\zeta_{ k})\|_{\L_{q_{k+1},p_{k+1}} }+
\| u \|_{\L_{q_{k+1},p_{k+1}}(C_{\rho_{k+1}})}\big).
\end{equation}
By Theorem \ref{theorem 4.3,1} we can 
continue as
$$
\|u\|_{\L_{q_{k},p_{k}}(C_{\rho_{k}})}
\leq N\big(\|\cL_{0}(u\zeta_{ k})\|_{\L_{q_{k+1},p_{k+1}} }+
\| u \|_{\L_{q_{k+1},p_{k+1}}(C_{\rho_{k+1}})}\big)
$$
$$
\leq N\|Du|D\zeta_{k}|\,\|_{\L_{q_{k+1},p_{k+1}}(C_{\rho_{k+1}})}+
\| u \|_{\L_{q_{k+1},p_{k+1}}(C_{\rho_{k+1}})}\big),
$$
where we used that $\cL_{0} u=0$. In  light of Remark \ref{remark 8.12.1}
$$
\|Du|D\zeta_{k}|\,\|_{\L_{q_{k+1},p_{k+1}}(C_{\rho_{k+1}})}\leq N
\|Du \|_{\L_{q_{k+1},p_{k+1}}(C_{\rho_{k+1/2}})}
$$
$$
\leq N
\|u \|_{\L_{q_{k+1},p_{k+1}}(C_{\rho_{k+1}})}.
$$
Hence, for $k
=0,1,...,n-1$
$$
\|u\|_{\L_{q_{k},p_{k}}(C_{\rho_{k}})}
\leq N\| u \|_{\L_{q_{k+1},p_{k+1}}(C_{\rho_{k+1}})},
$$
and \eqref{8.11.2} follows. \qed

\begin{theorem}
                    \label{theorem 4.29,4}
Let $p,q\in (1,\infty),\beta\in (1,\infty) $,
$u\in \dot \E^{1,2}_{q,p, \beta}$. 
Suppose that condition \eqref{4.3,5}
is satisfied with $\hat a_{0} $ from \eqref{8.22.4}
when $r=p\wedge q$.
Then
\begin{equation}
                        \label{4.29,3}
\|\partial_{t}u,D^{2}u\|_{\dot \E_{q,p,\beta}}
\leq N(d,\delta,p,q,\beta)\|f\|_{\dot \E_{q,p,\beta}},
\end{equation}
where $f=\cL_{0} u$.
\end{theorem}

Proof. First we observe that it suffices to
prove \eqref{4.29,3} for $u$ with compact support. For such $u$ and   $\lambda\in(0,1]$
define
$$
f_{\lambda}=\cL_{0} u-\lambda u,\quad g_{\lambda}=f_{\lambda}I_{C_{4}},\quad h_{\lambda}=f_{\lambda}I_{C_{4}^{c}}
$$ 
and let $G_{\lambda},H_{\lambda}$ be solutions
of class $\W^{1,2}_{q,p}(\bR^{d+1})$
of the equations 
\begin{equation}
                        \label{4.29,30}
\cL_{0} G_{\lambda}-\lambda G_{\lambda}=g_{\lambda},\quad 
\cL_{0} H_{\lambda}-\lambda H_{\lambda}=h_{\lambda}. 
\end{equation}
Since $f_{\lambda}$ is compactly supported
element of $\dot \E_{q,p,\beta}$,
it belongs to $\L_{q,p}$ and $G_{\lambda},
H_{\lambda}$ exist by Theorem 
\ref{theorem 4.3,3}. Furthermore,
$f_{\lambda}\in \dot E_{p\wedge q,\beta}
(\bR^{d+1})$,
so that by Theorem \ref{theorem 4.30,1}
the equations in \eqref{4.29,30} have
solutions in $\dot E^{1,2}_{p\wedge q,
 \beta }(\bR^{d+1})$.
These are the same $G_{\lambda},H_{\lambda}$, since these solutions
in different spaces were obtained
by the method of continuity
starting from the heat equation
for which the explicit formulas
for solutions are available.

Next, note that $\hat H_{\lambda}
:=e^{t}H_{\lambda}$ satisfies
$$
\cL _{0}\hat H_{\lambda}-(\lambda+1)
\hat H_{\lambda}=e^{t}h_{\lambda}:=\hat
h_{\lambda}
$$
and $\hat h_{\lambda}\in \L_{q,p}(\bR^{d+1}_{0})$ with its norm dominated by a constant independent of $\lambda$, since
$h_{\lambda}$ has compact support.
It follows from Theorem \ref{theorem 4.3,3}
that 
\begin{equation}
                  \label{1.5,2}
\sup_{\lambda\in(0,1]}
\|H_{\lambda}\|_{\L_{q,p}(\bR^{d+1}_{0})}
\leq \sup_{\lambda\in(0,1]}\|\hat H_{\lambda}\|_{\L_{q,p}(\bR^{d+1}_{0})}<\infty. 
\end{equation}

Concerning $G_{\lambda}$ note that
by Theorem \ref{theorem 4.3,3} 
$$
\| D^{2}G_{\lambda}\|_{\L_{q,p}(C_{1})}\leq N\|g_{\lambda}
\|_{\L_{q,p}}\leq N\|f_{\lambda}\|_{\L_{q,p}(C_{4})} \leq N\|f_{\lambda}\|_{\dot \E_{q,p,\beta}}. 
$$

Now observe that $\cL_{0}(e^{-\lambda t}H_{\lambda})=e^{-\lambda t}h_{\lambda}=0$ 
in $C_{4}$ and by Theorem \ref{theorem 8.11.1}
we have 
$$
\|e^{-\lambda t}H_{\lambda}-c\|_{L_{q,p}(C_{2})}
\leq N\|e^{-\lambda t}H_{\lambda}
-c\|_{L_{q\wedge p}(C_{4})}
$$
$$
\leq N\|H_{\lambda}
-c\|_{L_{q\wedge p}(C_{4})}+(1-e^{-16\lambda })
\|H_{\lambda}
 \|_{L_{q\wedge p}(C_{4})}
$$ 
$$
\leq N\|H_{\lambda}
-c\|_{L_{q\wedge p}(C_{4})}+N(1-e^{-16\lambda }),
$$ 
where $c$ is any affine function
and we used \eqref{1.5,2}. It follows by Theorem ~\ref{theorem 4.29,3} that
$$
\|D^{2}H_{\lambda}\|_{L_{q,p}(C_{1})}
\leq N\|e^{-\lambda t}H_{\lambda}-c\|_{L_{q,p}(C_{2})}\leq N\|H_{\lambda}
-c\|_{L_{q\wedge p}(C_{4})}+N(1-e^{-16\lambda }).
$$

 By Poincar\'e's
inequality for appropriate $c$
$$
\|H_{\lambda}
-c\|_{L_{q\wedge p}(C_{4})}\leq
N
\|\partial_{t}H_{\lambda},D^{2}H_{\lambda}\|_{L_{q\wedge p}(C_{4})}
\leq N
\|\partial_{t}H_{\lambda},D^{2}H_{\lambda}\|_{\dot E_{p\wedge q,\beta}}
$$
which, in turn, by Theorem
\ref{theorem 4.30,1} is dominated by 
$$
N\|h_{\lambda}\|_{\dot E_{p\wedge q,\beta}}\leq N\|f_{\lambda}\|_{\dot E_{p\wedge q,\beta}}\leq N\|f_{\lambda}\|_{\dot \E_{q,p,\beta}}.
$$

Since $u=G_{\lambda}+H_{\lambda}$
and \eqref{1.5,2} holds
it follows that 
$$
\|D^{2}u\|_{\L_{q,p}(C_{1})}
\leq \lim_{\lambda\downarrow 0}(
\|D^{2}G_{\lambda}\|_{\L_{q,p}(C_{1})}
+\|D^{2}H_{\lambda}\|_{\L_{q,p}(C_{1})}) 
$$
$$
\leq N\lim_{\lambda\downarrow 0}\|f_{\lambda}\|_{\dot \E_{q,p,\beta}}
\leq N\|f \|_{\dot \E_{q,p,\beta}}. 
$$
Parabolic scaling (preserving 
 \eqref{4.3,5})
shows that for any
$\rho>0$ 
$$
\rho^{\beta}\|D^{2}u\|_{\L_{q,p}(C_{\rho})}
\leq N\|f \|_{\dot \E_{q,p,\beta}} 
$$
and then shifting the origin leads to
\eqref{4.29,3} in what concerns $D^{2}u$.
Since $\partial_{t}u=f-a^{ij}D_{ij}u$,
the theorem is proved. \qed

\def\bim{\protect\eqref{4.30,2},
\protect\eqref{4.3,5}\,}
\mysection{Localizing  \bim and the proof
of Theorem \ref{theorem 5.8,20}}
                \label{subsection 8.22.1}
 
Here we argue under condition
\eqref{3.21.01}.

\begin{remark}
                \label{remark 4.14,1}
If $f(t,x)\geq 0$ is nonnegative, $r\in(0,\infty),\mu
\in(0,1)$, then
with 
$\zeta=|C_{r\mu}|^{-1}I_{C_{r\mu}}$ we have
$$
\dashnorm I_{C_{r}}f\|_{\L_{q,p}}
=N(d)r^{-d/p-2/q}\Big\|\int_{C_{(1+\mu)r}} 
  \zeta(s-\cdot,y-\cdot)I_{C_{r}}f(\cdot,\cdot)\,dyds\Big\|_{\L_{q,p}} 
$$
$$
\leq Nr^{-d/p-2/q}\int_{C_{(1+\mu)r}}\|\zeta(s-\cdot,y-\cdot)f\|_{\L_{q,p}}\,dyds   
$$
$$
\leq N(d)(1+\mu)^{d+2}   \mu  ^{d/p+2/q-(d+2)}\sup_{C\in
\bC_{\mu r}}\dashnorm
f\|_{\L_{q,p}(C )}. 
$$

Since $\mu\leq 1$ and $p,q>1$
it follows that for $\rho\leq r$ 
$$
 \sup_{C\in\bC_{r}}\dashnorm f\|_{\L_{q,p}(C)}\leq N(d) (\rho/r)^{d/p+2/q-(d+2)} \sup_{C\in\bC_{\rho}}\dashnorm f\|_{\L_{q,p}(C)}. 
$$
\end{remark}

\begin{lemma}
                          \label{lemma 5.7,1}
There exist a   constant $\kappa
=\kappa(d )\in(0,1)$ and a  
constant $K=K(d )$ such that for any $\rho>0$  there
exists an $\bS_{\delta}$-valued $\check a$
on $\bR^{d+1}$
such that $a=\check a$ in $C_{\kappa \rho }$
and $\check a^{ \shharp}\leq 
Ka_{\rho}^{ \shharp}$.
\end{lemma}
  
This lemma is proved in \cite{Kr_25_1}.
 
Recall that we fixed
$$
\rho_{a},\rho_{b}\in(0,1].  
$$

Let $\zeta_{0}\in C^{\infty}_{0}$
have suport in $C_{\kappa \rho_{a}} $ and
equal $1$ in $(1/2)C_{\kappa \rho_{a}}$
(with dimensions $(1/2)\kappa \rho_{a}$ in
the $x$-space and $(1/4)\kappa^{2} \rho^{2}_{a}$
in time and the same center),
where $\kappa$ is taken from Lemma \ref{lemma 5.7,1}. 
\begin{lemma}
                     \label{lemma 5.8,1}
There exists
$\hat a=\hat a(d,\delta,p,q,\beta)>0$
such that, if $a^{\shharp}_{\rho_{a}}\leq 
\hat  a$, 
then for any $u\in \E^{1,2}_{q,p,\beta}$ 
\begin{equation}
                        \label{5.8.3}
\|\zeta_{0}\partial_{t}u,\zeta_{0}D^{2}u\|_{  \E_{q,p,\beta}}
\leq N  \|f\|_{  \E_{q,p,\beta}}+ N_{1}\|Du,u\|_{  \E_{q,p,\beta}} , 
\end{equation}
where $f=\cL_{0} u$, $N =N (d,\delta,p,q,\beta )$,
$N_{1} =N (d,\delta,p,q,\beta,\rho_{a} )$. 
\end{lemma}

Proof. Observe that
\begin{equation}
                        \label{5.8,2}
\|\zeta_{0}\partial_{t}u,\zeta_{0}D^{2}u\|
_{  \E_{q,p,\beta}}\leq
\|\partial_{t}(\zeta_{0}u),D^{2}(\zeta_{0}u)\|_{\dot
\E_{q,p,\beta}}+N\|Du,u\|_{ \E_{q,p,\beta}}
\end{equation}
and, owing to Lemma \ref{lemma 5.7,1}
and Theorem \ref{theorem 4.29,4}, the first term is estimated through
$$
(\partial_{t}+\check a^{ij}D_{ij})(\zeta_{0}u)
=\cL_{0}(\zeta_{0}u)=\zeta_{0}f + g,
$$
where, thanks to  $\beta\leq d/p+2/q$, the $\dot \E_{q,p,\beta}$-norm of $g$
is estimated by the last term in
\eqref{5.8,2}. \qed

Naturally, \eqref{5.8.3} holds for
any shift of $\zeta_{0}$ and this
along with Remark \ref{remark 4.14,1} proves
the following with    $ \|Du,u\|_{  \E_{q,p,\beta}}$
in place of $ 
\|u\|_{  \E_{q,p,\beta}}$. The term 
$\|Du \|_{  \E_{q,p,\beta}}$ is eliminated
by using interpolation inequalities.

\begin{theorem}
                     \label{theorem 5.8,1}
There exists
$\hat a=\hat a(d,\delta,p,q,\beta)>0$
such that, if $a^{\shharp}_{\rho_{a}}\leq 
\hat  a$, then for any $u\in \E^{1,2}_{q,p,\beta}$ 
\begin{equation}
                        \label{5.8.40}
\| \partial_{t}u, D^{2}u\|_{  \E_{q,p,\beta}}
\leq N  \|\cL_{0} u\|_{  \E_{q,p,\beta}}+N\| u\|_{  \E_{q,p,\beta}} , 
\end{equation}
where $N =N (d,\delta,p,q,\beta,\rho_{a})$.
 
\end{theorem}

Next, as a few times before, we use Agmon's
idea along
with the interpolation inequalities to treat $\cL_{0}-\lambda$. 
After that to deal with the first-order
terms, note that 
for $\zeta\in C^{\infty}_{0}(
2C_{1})$ such that $\zeta=1$ on $C_{1}$,
by H\"older's inequality
 we obtain 
for $\rho\leq 1$ and $u\in C^{\infty}_{0}$ that 
$$
\rho^{\beta}\|b^{i} D_{i}u\|_{\L_{q,p}(C_{\rho})}
\leq  \|\zeta b^{i} D_{i}(\zeta u)\|_{\dot \E_{q,p,\beta}(C_{\rho})}
\leq N\|b  \|_{ \E_{q\beta,p\beta,1}}\|   D(\zeta u)  \|_{\dot \E_{r_{1},r_{2},\beta-1}},
$$
where $r_{1}(\beta-1)=q\beta,r_{2}(\beta-1)=p\beta$. 
By the embedding theorem (see Corollary 5.7 of \cite{Kr_23} in which the order of integration
in the mixed-norm space did not play any
crucial role) $\|   D(\zeta u)  \|_{\dot \E_{r_{1},r_{2},\beta-1}}\leq N\|\zeta u\|_{ \dot\E^{1,2}_{q,p,\beta }}$,
implying that
\begin{equation}
                             \label{5.8,5}
 \|b^{i} Du\|_{\E_{q,p,\beta}}
\leq N\|b  \|_{ \E_{q\beta,p\beta,1}}
\|u\|_{ \E^{1,2}_{q,p,\beta }}. 
\end{equation}
 
By considering the cut-off mollifiers of
$u\in \E^{1,2}_{q,p}$ we extend 
\eqref{5.8,5} to all $u\in \E^{1,2}_{q,p}$.
In this way we come to the a priori estimate
below.  Note 1 in the condition on $\bar b$.

\begin{theorem}
                     \label{theorem 5.8,2}
There 
exist
$$
\hat a=\hat a(d,\delta,p,q,\beta)>0,\quad\hat b=\hat b(d,\delta,p,q,\beta,\rho_{a} )>0,
$$ 
$$\lambda_{0}= \hat \lambda(d,\delta,p,q,\beta,\rho_{a} ) >0,\quad
N_{0}=N_{0}(d,\delta,p,q,\beta,\rho_{a}  ),
$$
such that, if 
\begin{equation}
                              \label{8.22.6}
a^{\shharp}_{\rho_{a}}\leq 
\hat  a,\quad \bar b_{q\beta,p\beta,1} \leq \hat b,
\end{equation}
 then for any $u\in \E^{1,2}_{q,p,\beta}$  
 and $\lambda\geq\lambda_{0}$ 
\begin{equation}
                        \label{5.8.4}
\|\lambda u,\sqrt\lambda Du, \partial_{t}u, D^{2}u\|_{  \E_{q,p,\beta}}
\leq N_{0} \|f\|_{  \E_{q,p,\beta}}, 
\end{equation}
where   $f=\cL_{0} u+b^{i}D_{i}u-\lambda u$.
Furthermore, for any
$f\in \E_{q,p,\beta}$ and $\lambda\geq\lambda_{0}$ there exists a unique
$u\in \E^{1,2}_{q,p,\beta}$ such that in
$\bR^{d+1}$
$$
\cL_{0} u+b^{i}D_{i}u-\lambda u=f.
$$
 
\end{theorem}

Proof. Of course, in condition \eqref{8.22.6} we take $\hat a$
from Theorem \ref{theorem 5.8,1} and find
$\hat b$, $\lambda_{0}$, $N_{0}$ as described
before the statement of the theorem.
To prove the existence of solutions we use
the method of continuity. In light of the a priori estimate  \eqref{5.8.4},
for this method to work, we  need only prove
the solvability for the heat equation
$$
\partial_{t}u+\Delta u-\lambda u=f.
$$
If $f$ is bounded and with each derivative  bounded, the solution $u$
admits an explicit representation,
is bounded and has bounded
derivatives.
Therefore, $u\in \E^{1,2}_{q,p,\beta}$
and admits estimate \eqref{5.8.4}.
For general $f\in \E_{q,p,\beta}$
it suffices to use its mollifiers.  \qed

{\bf Proof of Theorem \ref{theorem 5.8,20}}. 
We take the same $\hat a$
$\hat b$, $\lambda_{0}$, $N_{0}$ as in the above proof.
Recall that $\rho_{b}\in(0,1]$. Owing to Theorem \ref{theorem 5.8,2}
our assertions are true if $\rho_{b}=1$.
Then, having \eqref{5.10.2} for $\rho_{b}=1$,
for $\rho_{b}<1$ we make a parabolic dilation
(not affecting $a^{\shharp}_{\rho_{a}}\leq 
\hat  a$) and see that for $\lambda\geq\lambda_{0}\rho_{b}^{-2}$, where
$\lambda_{0}$ is from Theorem~\ref{theorem 5.8,2},
$$
\sup_{\rho\leq\rho_{b}}\rho^{\beta}\sup
_{C\in\bC_{\rho}}\dashnorm\lambda u,\sqrt\lambda Du, \partial_{t}u, D^{2}u\|_{ \L_{q,p}(C)} 
$$
\begin{equation}
                        \label{5.8.5}
\leq N_{0}\sup_{\rho\leq\rho_{b}}\rho^{\beta}\sup
_{C\in\bC_{\rho}}\dashnorm  f\|_{ \L_{q,p}(C)}
\leq N_{0}\|f\|_{  \E_{q,p,\beta}}. 
\end{equation}
By Remark \ref{remark 4.14,1} the left-hand
side of \eqref{5.8.5} times $N(d)\rho_{b}^{-\alpha}$ dominates the left-hand side
of \eqref{5.10.2}. This proves the theorem.
\qed

\end{document}